\documentclass[10pt,reqno]{amsart}
\usepackage{amscd,amsmath,amsthm,amsfonts,latexsym,amssymb}

\theoremstyle{plain}
\newtheorem{theo+}           {Theorem}      [section]
\newtheorem{prop+}  [theo+]  {Proposition}
\newtheorem{coro+}  [theo+]  {Corollary}
\newtheorem{lemm+}  [theo+]  {Lemma}

\theoremstyle{definition}
\newtheorem{rema+}  [theo+]  {Remark}
\newtheorem{defi+}  [theo+]  {Definition}

\newenvironment{theorem}{\begin{theo+}}{\end{theo+}}
\newenvironment{proposition}{\begin{prop+}}{\end{prop+}}
\newenvironment{corollary}{\begin{coro+}}{\end{coro+}}
\newenvironment{lemma}{\begin{lemm+}}{\end{lemm+}}
\newenvironment{remark}{\begin{rema+}}{\end{rema+}}

\begin{document}
\numberwithin{equation}{section}

\title[Annihilating fields of standard modules]{Annihilating fields
  of standard modules \\ for affine Lie algebras}
\author[Julius Borcea]{Julius Borcea}
\address{Institut de Recherche Math\'ematique Avanc\'ee, Universit\'e Louis Pasteur, 7 rue Ren\'e Descartes, 67084 Strasbourg Cedex, France}
\email{borcea@math.u-strasbg.fr}

\begin{abstract}
  Given an affine 
  Kac-Moody Lie algebra $\tilde{\mathfrak{g}}[\sigma]$ of arbitrary type, we 
  determine certain minimal 
  sets of annihilating fields of standard
  $\tilde{\mathfrak{g}}[\sigma]$-modules. 
  We then use these sets in order to obtain a 
  characterization of standard $\tilde{\mathfrak{g}}[\sigma]$-modules
  in terms of irreducible loop $\tilde{\mathfrak{g}}[\sigma]$-modules,
  which proves to be a useful tool for combinatorial constructions of bases 
  for standard $\tilde{\mathfrak{g}}[\sigma]$-modules.
\end{abstract}

\maketitle

\section*{Introduction}

It is a well-known fact that 
suitably specialized characters of standard modules
for affine Lie algebras can be expressed as certain infinite
products. In turn, these products may be interpreted as generating functions 
of partition functions for colored 
partitions defined by congruence conditions (cf.\@~\cite{An}). Such arguments 
have revealed the connections existing between the character theory of 
standard modules for affine Lie algebras and (the product side of) partition 
identities of Rogers-Ramanujan type (cf.\@~\cite{LM}). It was therefore natural 
to ask whether these identities could be understood through representation 
theory as manifestations of a new algebraic structure that would lead to bases
  (of the corresponding modules) which would be parametrized by (colored 
partitions satisfying) difference conditions. This question has been answered 
in the affirmative in a range of cases -- including the original 
Rogers-Ramanujan identities -- by Lepowsky and Wilson (\cite{LW1,LW2}). In 
order to obtain this result, they introduced and studied a new kind of 
associative algebras, the so-called $\mathcal{Z}$-algebras (see also 
\cite{LP}). (It turned out later on that these structures were closely 
related to the parafermion algebras that appeared in physics literature in 
the mid-eighties (cf.\@~\cite{DL}).) Further applications of the 
$\mathcal{Z}$-algebra theory to the problem of constructing (what came to be 
called) fermionic representations may be found in e.\@~g. 
\cite{Ca2,Ca3,MP1,P2,TX}.  

In this paper we make use of vertex operator techniques in order to determine 
certain sets of annihilating fields of standard $\tilde{\mathfrak{g}}[\sigma]$
-modules for an affine Kac-Moody Lie algebra $\tilde{\mathfrak{g}}[\sigma]$ of 
arbitrary type. It turns out that these 
sets have a structure of loop $\tilde{\mathfrak{g}}
[\sigma]$-modules of level 0 (see \cite{CP} for the definition). Moreover, 
these sets are 
minimal in the sense that they are irreducible loop $\tilde{\mathfrak{g}}
[\sigma]$-modules (Theorems 2.9 and 2.1). We also get a description of 
the maximal submodule of the Verma $\tilde{\mathfrak{g}}[\sigma]$-module
$M(\Lambda)$, $\Lambda \in P_{+}$, in terms of these loop
$\tilde{\mathfrak{g}}[\sigma]$-modules (Theorem 2.13). As a
consequence, we obtain a characterization of standard
$\tilde{\mathfrak{g}}[\sigma]$-modules by means of irreducible loop
$\tilde{\mathfrak{g}}[\sigma]$-modules (Theorem 2.14). Partial analogues of 
this result have proven to be very useful 
for combinatorial constructions of bases for standard modules in several 
particular cases (cf., e.\@~g., \cite{Ca1,MP2,MP3}). 
The main ingredients used below are a combination between the untwisted  
representation theory of affine VOAs (\cite{DL,Li1,MP2})
and Kac's classification of finite-order automorphisms of
finite-dimensional simple Lie algebras (\cite{K}), together with Li's 
results on deformations of vertex operator maps by means of certain 
endomorphisms associated with inner automorphisms (\cite{Li2}).

The results presented here generalize those obtained in the untwisted 
case in \cite[\S 5]{MP2} and \cite{P1}, and they may be of 
interest in their own right. Moreover, they may also provide an appropriate 
setting for the study of standard modules for affine Lie algebras of arbitrary 
type. Indeed, the above-mentioned results led in the untwisted case to the 
construction of bases for standard modules of arbitrary level 
for the rank two untwisted affine Lie algebra $A_{1}^{(1)}$ in 
\cite{MP2} (bases for these same modules had been previously constructed in 
\cite{LP}). As a consequence, a series of combinatorial identities of 
Rogers-Ramanujan type was obtained (see also \cite{MP3} for further 
applications). The representation theories of the rank two 
affine Lie algebras $A_{1}^{(1)}$ and $A_{2}^{(2)}$ are to a certain extent 
prototypical for the representation theory of all the untwisted respectively 
twisted 
affine Lie algebras. It seems therefore natural to investigate whether, for 
instance, standard $A_{2}^{(2)}$-modules of arbitrary level could be dealt 
with in a similar fashion by using the setting developed here. If successful, 
this approach may lead in particular to a combinatorial explanation of the 
duality-like property for rank two affine Lie algebras discussed in \cite{Bo}  
and is currently under study.

It is worth
mentioning that a qualitative version of the results in \cite[\S 5]{MP2} 
was recently obtained in \cite{FM} for the class of admissible 
representations (cf.\@~\cite{KW}) of untwisted affine Lie algebras. We believe
that our results can be extended to admissible 
representations of affine Lie algebras of arbitrary 
type, which may then lead to further applications (like those in e.\@~g. 
\cite{Ad}). 

\smallskip

\noindent
{\bf Acknowledgements.} I would like to thank Mirko Primc and Arne Meurman 
for numerous stimulating discussions, and the anonymous referee for his 
careful reading of this paper.

\section{Preliminaries and notations}

\subsection{Realizations of affine Lie algebras.} 

Let $\mathfrak{g}$ be a finite-dimensional simple complex Lie algebra. 
Fix a $CSA$ $\mathfrak{t}$ of $\mathfrak{g}$ and let $\mu$ be an 
automorphism of $\mathfrak{g}$ of order $r$ ($=1,2,$ or 3) induced by
an automorphism of order $r$ of the Dynkin diagram of $\mathfrak{g}$
with respect to $\mathfrak{t}$. Let
$\varepsilon$ be a primitive $r$-th root of unity and denote by
$\mathfrak{g}_{[i]}$ the $\varepsilon^{i}$-eigenspace of $\mu$ in
$\mathfrak{g}$, $i\in \mathbf{Z}_{r}$. Then 
$\mathfrak{g}_{[0]}$ is a simple subalgebra of
$\mathfrak{g}$, the space $\mathfrak{t}_{[0]}:=\mathfrak{g}_{[0]}\cap 
\mathfrak{t}$ is
a $CSA$ of $\mathfrak{g}_{[0]}$, and the $\mathfrak{g}_{[0]}$-modules
$\mathfrak{g}_{[1]}$ and $\mathfrak{g}_{[-1]}$ are irreducible and
contragredient. Set $l=\text{rank}\,\mathfrak{g}_{[0]}$, and let 
$\{\beta_{1},\ldots,\beta_{l}\}\subset 
\mathfrak{t}_{[0]}^{*}$ be a basis of the root system of 
$\mathfrak{g}_{[0]}$ and $\{E_{j}, F_{j}, H_{j}\mid j\in
\{1,\ldots,l\}\}$ a corresponding set of
canonical generators of $\mathfrak{g}_{[0]}$. Let $\beta_{0}\in
\mathfrak{t}_{[0]}^{*}$ be the lowest weight of the $\mathfrak{g}_{[0]}$-module
$\mathfrak{g}_{[1]}$, and let $E_{0}$ and $F_{0}$
be a lowest weight vector of the $\mathfrak{g}_{[0]}$-module
$\mathfrak{g}_{[1]}$ respectively a highest weight vector of the 
$\mathfrak{g}_{[0]}$-module
$\mathfrak{g}_{[-1]}$. We assume that $E_{0}$ and $F_{0}$ are normalized so 
that $[H_{0},E_{0}]=2E_{0}$,
where $H_{0}=[E_{0},F_{0}]$. For $i,j\in \{0,1,\ldots,l\}$ define
$a_{ij}=\beta_{j}(H_{i})$ and let $A=(a_{ij})_{i,j=0}^{l}$. Then $A$ is a 
GCM of affine type, and so there are 
positive integers $a_{0},\ldots,a_{l}$ such 
that $(a_{0},\ldots,a_{l})A^{t}=0$. 
Equivalently, there exist positive integers
$a\check{}_{0},\ldots,a\check{}_{l}$ such that
$(a\check{}_{0},\ldots,a\check{}_{l})A=0$. 
Both these sets of integers are assumed to be normalized so that 
$\gcd(a_{0},\ldots,a_{l})=\gcd(a\check{}_{0},\ldots,a\check{}_{l})$ $=1$. Then 
$h:=\sum_{j=0}^{l}a_{j}$ and
$h\check{}:=\sum_{j=0}^{l}a\check{}_{j}$ are the Coxeter
number respectively the dual Coxeter
number of the matrix $A$, which will be denoted by $X_{N}^{(r)}$ if 
$\mathfrak{g}$ is of type $X_{N}$
($X=A,B,\ldots,G$ and $N\ge 1$). We shall use the Dynkin diagrams of 
the affine GCMs as
listed in \cite{KKLW}, that is, with the vertex corresponding to the 0-th
index always occurring at the left-end of the diagram. 

Let $\mathbf{s}=(s_{0},s_{1},\ldots,s_{l})$ be a sequence of
nonnegative relatively prime integers, and set
$T=r\sum_{j=0}^{l}s_{j}a_{j}$. If $\eta$ is a primitive $T$-th root of
unity, the conditions $\nu(H_{j})=H_{j}$, $\nu(E_{j})=\eta^{s_{j}}E_{j}$, 
$0\le j\le l$, define a $T$-th order automorphism $\nu$ of
$\mathfrak{g}$, the so-called $\mathbf{s}$-automorphism. Since $\mathbf{C}$ 
is algebraically
closed, every finite-order automorphism of
$\mathfrak{g}$ arises in this way up to conjugation by an automorphism
of $\mathfrak{g}$ (\cite[Theorem 8.6]{K}). Notice that
in this terminology, the original diagram automorphism $\mu$ becomes
the $(1,0,\ldots,0)$-automorphism of $\mathfrak{g}$. Denote
by $\mathfrak{g}_{(j)}$ the $\eta^{j}$-eigenspace of $\nu$ in
$\mathfrak{g}$ for $j\in \mathbf{Z}_{T}$. The
$\mathbf{Z}_{T}$-gradation $\mathfrak{g}=\coprod_{j\in
  \mathbf{Z}_{T}}\mathfrak{g}_{(j)}$ is accordingly called the
$\nu$-gradation (or $\mathbf{s}$-gradation) and a graded subspace of 
$\mathfrak{g}$ is said
to be $\nu$-graded. Let $\langle \cdot\, ,\cdot \rangle$ be a nondegenerate
symmetric $\mathfrak{g}$-invariant bilinear form on
$\mathfrak{g}$. Being a multiple of the Killing form,
$\langle \cdot\, ,\cdot \rangle$ is also $\nu$-invariant and remains 
nonsingular on
the $CSA$ $\mathfrak{t}_{[0]}$ of $\mathfrak{g}_{[0]}$. We may
therefore identify $\mathfrak{t}_{[0]}$ with $\mathfrak{t}_{[0]}^{*}$ by 
means of the
restricted form. Furthermore, we may assume that $\langle \cdot\, ,\cdot 
\rangle$
is normalized so that $\langle \beta_{0},\beta_{0} \rangle =2a\check{}_{0}/r$,
which then implies that
$a\check{}_{j}=r\langle \beta_{j},\beta_{j} \rangle a_{j}/2$ for 
$j=0,1,\ldots,l$
(cf.\@~\cite[Proposition 1.1]{KKLW}). This normalization of the form
$\langle \cdot\, ,\cdot \rangle$ amounts to the condition that
$\langle \alpha,\alpha \rangle =2$ whenever $\alpha\in \mathfrak{t}^{*}$ 
is a long
root of $\mathfrak{g}$, in which case
the Killing form  equals $2h\check{}\langle \cdot\, ,\cdot \rangle$. Define 
the Lie algebras
\begin{equation}\label{1.1}   
  \hat{\mathfrak{g}}[\nu]=\oplus_{j=0}^{T-1}\mathfrak{g}_{(j)}\otimes
  t^{\frac{j}{T}}\mathbf{C}[t,t^{-1}]\oplus \mathbf{C}c, \quad  
  \tilde{\mathfrak{g}}[\nu]= \hat{\mathfrak{g}}[\nu]\rtimes
  \mathbf{C}d,
\end{equation}
by the conditions
\begin{equation}\label{1.2}
  \begin{split}
  & c \;\, \mbox{central},\; c\neq 0, \; [d,a\otimes t^{m}]=m\, a\otimes
  t^{m},\\
  & [a\otimes t^{m},b\otimes t^{n}]=[a,b]\otimes
  t^{m+n}+m\delta_{m+n,0}\langle a,b \rangle c,
  \end{split}
\end{equation}
for $m,n\in \frac{1}{T}\mathbf{Z}$, $a\in \mathfrak{g}_{(mT \bmod T)}$,
$b\in \mathfrak{g}_{(nT \bmod T)}$. For $a\in \mathfrak{g}_{(j)}$, $n\in
\mathbf{Z}$, we shall frequently use $a(n+j/T)$ and
$\mathfrak{g}(n+j/T)$ to denote $a\otimes t^{n+j/T}$ and 
$\mathfrak{g}_{(j)}\otimes
  t^{n+j/T}$ respectively, and we often identify
  $\mathfrak{g}_{(0)}(0)$ with $\mathfrak{g}_{(0)}$. The space
  $\mathfrak{h}:=\mathfrak{t}_{[0]}\oplus \mathbf{C}c$ (res\-pectively
  $\mathfrak{h}^{e}:=\mathfrak{h}\rtimes \mathbf{C}d$) is a $CSA$ of 
$\hat{\mathfrak{g}}[\nu]$ (respectively $\tilde{\mathfrak{g}}[\nu]$).
 Let $\delta\in \mathfrak{h}^{e^{*}}$ be such that
$\delta|_{\mathfrak{h}}=0$, $\delta(d)=1$, and define 
$\alpha_{j}\in \mathfrak{h}^{e^{*}}$ by
$\alpha_{j}|_{\mathfrak{t}_{[0]}}=\beta_{j}$, $\alpha_{j}(c)=0$,
$\alpha_{j}(d)=s_{j}T^{-1}$ if $j=1,\ldots,l$, and
$\alpha_{0}=r^{-1}\delta-\sum_{j=1}^{l}a_{j}\alpha_{j}$, so that in
particular $\alpha_{0}(d)=s_{0}T^{-1}$. For $j\in \{0,1,\ldots,l\}$ define 
also 
\begin{equation}\label{1.3}
  e_{j}=E_{j}\otimes t^{\frac{s_{j}}{T}},\quad f_{j}=F_{j}\otimes
  t^{-\frac{s_{j}}{T}},\quad
  h_{j}=H_{j}+\frac{2s_{j}}{T\langle \beta_{j},\beta_{j} \rangle}c.
\end{equation}
Then $\{e_{j},f_{j},h_{j},d\mid j\in \{0,1,\ldots,l\}\}$ is a system
of canonical generators of $\tilde{\mathfrak{g}}[\nu]$, viewed as the
$\nu$-twisted affine Kac-Moody Lie algebra of rank $l+1$ associated to the GCM 
$A$, in what is called the 
$\mathbf{s}$-realization of this algebra. Note that the canonical central 
element $c$ equals $\sum_{j=0}^{l}a\check{}_{j}h_{j}$, and let
$\tilde{\mathfrak{g}}[\nu]_{i}:=\{a\in \tilde{\mathfrak{g}}[\nu]\mid
[d,a]=i\,a\}$, $i\in \frac{1}{T}\mathbf{Z}$. The corresponding 
$\frac{1}{T}\mathbf{Z}$-gradation
$\tilde{\mathfrak{g}}[\nu]=\coprod_{i\in
  \frac{1}{T}\mathbf{Z}}\tilde{\mathfrak{g}}[\nu]_{i}$ is then called the
$\nu$-gradation of $\tilde{\mathfrak{g}}[\nu]$. 
Let as usual $\mathfrak{n}_{+}$ and $\mathfrak{n}_{-}$ denote the
subalgebras of $\hat{\mathfrak{g}}[\nu]$ generated by
$e_{0},\ldots,e_{l}$ and by $f_{0},\ldots,f_{l}$ respectively, so that 
one has the triangular decompositions
\begin{equation}\label{1.4}
  \hat{\mathfrak{g}}[\nu]=\mathfrak{n}_{-}\oplus \mathfrak{h}\oplus
  \mathfrak{n}_{+},\quad \tilde{\mathfrak{g}}[\nu]=\mathfrak{n}_{-}\oplus 
\mathfrak{h}^{e}\oplus
  \mathfrak{n}_{+},
\end{equation}
and corresponding decompositions of the universal enveloping algebras of 
$\hat{\mathfrak{g}}[\nu]$ and $\tilde{\mathfrak{g}}[\nu]$.

A $\hat{\mathfrak{g}}[\nu]$-module $V$ is said to be restricted 
and 
of level $l$ if $\mathfrak{g}_{(j)}(n+j/T)\cdot v=0$ for any
$v\in V$ and $n\gg 0$ and $c$ acts as $l\mbox{id}_{_{V}}$ on
$V$. In particular, any highest-weight module is restricted. 
Given a restricted $\hat{\mathfrak{g}}[\nu]$-module $V$ and $a\in
\mathfrak{g}_{(j)}$, we shall consider the generating function of
operators on $V$
\begin{equation}\label{1.5}
  a(\nu;z)=\sum_{n\in \mathbf{Z}}a(n+j/T)z^{-n-\frac{j}{T}-1}\in
  (\mbox{End}\,V)\big[\big[z^{1/T},z^{-1/T}\big]\big].
\end{equation}
We shall sometimes write $a_{n+j/T}$ when we think of $a(n+j/T)$ as a
coefficient of $a(\nu;z)$, and $a(\mbox{id}_{\mathfrak{g}};z)$ will be
denoted simply by $a(z)$. Let $M(\Lambda)$ be the Verma 
$\tilde{\mathfrak{g}}[\nu]$-module with highest weight 
$\Lambda\in \mathfrak{h}^{e^{*}}$. Denote by $M^{1}(\Lambda)$ its unique 
maximal proper submodule and let $L(\Lambda)=M(\Lambda)/M^{1}(\Lambda)$. 
Recall that a highest-weight $\tilde{\mathfrak{g}}[\nu]$-module 
$V$ with highest weight $\Lambda$ is standard if there exists 
$m\ge 1$ such that $f_{i}^{m}\cdot v_{\Lambda}=0$, $i=0,1,\ldots,l$, where
$v_{\Lambda}\in V$ is a highest weight vector, which then implies that
$\Lambda\in P_{+}:=\{\Lambda\in \mathfrak{h}^{e^{*}}\mid \Lambda(h_{i})\in
  \mathbf{Z}_{_{\ge 0}}\; \mbox{for}\; i=0,\ldots,l\}$. Conversely, if 
$\Lambda\in P_{+}$ then 
$L(\Lambda)$ is a standard $\tilde{\mathfrak{g}}[\nu]$-module and by 
\cite[Corollary 10.4]{K} one has that
\begin{equation}\label{1.6}
  M^{1}(\Lambda)=\sum_{i=0}^{l}U(\mathfrak{n}_{-})f_{i}^{\Lambda(h_{i})+1}
\cdot v_{\Lambda}.
\end{equation}
Let finally $\Lambda_{i}\in P_{+}$, $0\le i\le l$, denote the fundamental 
weights determined by $\Lambda_{i}(h_{j})=\delta_{ij}$, $\Lambda_{i}(d)=0$.

\subsection{VOAs and modules.}

We refer to \cite{B,FHL,FLM} for the definition 
of a vertex (operator) algebra, and to \cite{DL,Li1,Li2} 
for the different notions of weak module for a VOA. The definition of a 
twisted VOA-module used below may be found in e.\@~g. ~\cite{Li2}. 

Let $(V,Y,\mathbf{1},\omega)$ be a VOA, and recall that $\mbox{id}_{_{V}}$ 
together with the component 
ope\-rators of the field $Y(\omega,z)=\sum_{n\in \mathbf{Z}}L(n)z^{-n-2}$ 
generate a representation of the 
Virasoro algebra on $V$. Let further $\sigma$ be an automorphism of order $T$ 
of $V$ and $V^{k}=\{a\in V\mid
  \sigma (a)=\exp (2k\pi i/T)a \}$, $0\le k\le T-1$, so that 
$V=\oplus_{k=0}^{T-1}V^{k}$. If $(M,Y_{M})$ is a $\sigma$-twisted $V$-module 
and $a\in V^{k}$, $b\in V$, then 
\begin{equation}\label{1.7}
Y_{M}(L(-1)b,z)=\frac{d}{dz}Y_{M}(b,z), 
\end{equation}
and one has the following consequences of the defining 
axioms (cf., e.\@~g., \cite{Li2}):
{\allowdisplaybreaks
  \begin{align}
    & z^{\frac{k}{T}}Y_{M}(a,z)\in (\mbox{End}\, M)[[z,z^{-1}]],\label{1.8}\\
    & [Y_{M}(a,z_{1}),Y_{M}(b,z_{2})]=\sum_{j=0}^{\infty}\frac{1}{j!}\!
\left(\!\left(\!\frac{\partial}{\partial z_{2}}\!\right)^{j}z_{1}^{-1}\delta
\!\left(\!\frac{z_{2}}{z_{1}}\!\right)\left(\!\frac{z_{2}}{z_{1}}\!\right)
^{\frac{k}{T}}\!\right)\!Y_{M}(a_{j}b,z_{2}),\label{1.9}  \\
    & Y_{M}(Y(a,z_{0})b,z_{2})=\mbox{Res}_{z_{1}}\left(\!\frac{z_{1}-z_{0}}
{z_{2}}\!\right)^{\frac{k}{T}}\!\Biggl[z_{0}^{-1}\delta\!\left(\!
\frac{z_{1}-z_{2}}{z_{0}}\!\right)\!Y_{M}(a,z_{1})Y_{M}(b,z_{2})\Biggr.\label{1.10}\\
    & \phantom{Y_{M}(Y(a,z_{0})b,z_{2})=\mbox{Res}_{z_{1}}\left(\!
\frac{z_{1}-z_{0}}{z_{2}}\!\right)^{\frac{k}{T}}}\quad 
    \Biggl.-z_{0}^{-1}\delta\!\left(\!\frac{-z_{2}+z_{1}}{z_{0}}\!\right)\!
Y_{M}(b,z_{2})Y_{M}(a,z_{1})\Biggr]. \nonumber 
\end{align}

We now describe briefly the so-called affine VOAs.} Let $\mathfrak{g}$ be a finite-dimensional simple Lie algebra with the form $\langle \cdot,\!\cdot
\rangle$ normalized as in \S 1.1, and form the untwisted
affine Kac-Moody algebra $\tilde{\mathfrak{g}}=\hat{\mathfrak{g}}\rtimes
  \mathbf{C}d$ 
as in \eqref{1.1}-\eqref{1.2}. Set 
  $\tilde{\mathfrak{g}}_{\ge 0}=\oplus_{n\ge 0}\mathfrak{g}(n)\oplus
  \mathbf{C}c \oplus \mathbf{C}d$ 
and let $-h\check{}\neq l\in \mathbf{C}$. Recall from \S 1.1 the 
fundamental weight $\Lambda_{0}\in P_{+}$ and define a
$\tilde{\mathfrak{g}}_{\ge 0}$-module structure on $\mathbf{C}$ by 
$c\cdot 1=l,\, d\cdot 1=0,\, \mathfrak{g}(n)\cdot 1=0 \text{ for }
  n\ge 0$. We may then form the Weyl module (or generalized Verma module) 
$N(l\Lambda_{0})=U(\tilde{\mathfrak{g}})\otimes_{U(\tilde{\mathfrak{g}}_{\ge 0})}
\mathbf{C}$, which is the so-called vacuum representation of level $l$ of 
$\tilde{\mathfrak{g}}$. Note that $N(l\Lambda_{0})$ is a restricted
$\tilde{\mathfrak{g}}$-module such that $N(l\Lambda_{0})\cong
U\big(\!\oplus_{_{n<0}}\mathfrak{g}(n)\big)$ as vector spaces, and that we 
may identify 
$\mathfrak{g}(-1)\otimes 1$ with $\mathfrak{g}$. Set 
$\mathbf{1}=1\otimes 1\in N(l\Lambda_{0})$  
and define the element 
$$\omega=\frac{1}{2(l+h\check{})}\sum_{j=1}^{\dim
    \mathfrak{g}}a^{j}(-1)^{2}\mathbf{1}\in N(l\Lambda_{0}),$$
where $\{a^{j}\mid j\in \{1,\ldots,\dim \mathfrak{g}\}\}$ is an orthonormal 
basis of $\mathfrak{g}$ with respect to $\langle \cdot,\! \cdot
\rangle$. Recall from \eqref{1.5} the series $a(z)$ in
this case and define the map 
\begin{gather*}
  Y:\, \mathfrak{g}(-1)\otimes 1 \longrightarrow (\mbox{End}\,
  N(l\Lambda_{0}))[[z,z^{-1}]]\\
  Y(a(-1)\otimes 1,z)=a(z), \text{ } a\in \mathfrak{g}.
\end{gather*}
One can show (cf., e.\@~g., \cite[Theorem 2.6]{MP2}) that
$Y$ extends uniquely to 
$N(l\Lambda_{0})$ in such a way that $N(l\Lambda_{0})$ becomes a VOA with
vacuum vector $\mathbf{1}$ and Virasoro 
element $\omega$ such that $\mathfrak{g}(-1)\otimes
1=N(l\Lambda_{0})_{1}$ (the weight one subspace of $N(l\Lambda_{0})$). 
Moreover, given any restricted $\hat{\mathfrak{g}}$-module $M$ of level
$l$, there is a canonical extension to $N(l\Lambda_{0})$ of the map
\begin{gather*}
  Y_{M}:\, \mathfrak{g}(-1)\otimes 1 \longrightarrow (\mbox{End}\,
  M)[[z,z^{-1}]]\\
  Y_{M}(a(-1)\otimes 1,z)=a(z), \text{ } a\in \mathfrak{g},     
\end{gather*}
such that $(M,Y_{M})$ becomes a weak $N(l\Lambda_{0})$-module
(\cite{DL,Li1,MP2}). Let finally $N^{1}(l\Lambda_{0})$ be
the unique maximal proper $\tilde{\mathfrak{g}}$-submodule of
$N(l\Lambda_{0})$ and notice that we may identify the irreducible
quotient $N(l\Lambda_{0})/N^{1}(l\Lambda_{0})$ with the 
$\tilde{\mathfrak{g}}$-module $L(l\Lambda_{0})$ defined in \S 1.1. We summarize
some of the results of the above-mentioned papers in
{\allowdisplaybreaks
\begin{theo+}
  For each $l\neq
  -h\check{}$, $(N(l\Lambda_{0}),Y,\mathbf{1},\omega)$ is a VOA of rank
  $\tfrac{l\dim \mathfrak{g}}{l+h\check{}}$ and any restricted 
$\hat{\mathfrak{g}}$-module of level
  $l$ is a weak $N(l\Lambda_{0})$-module. Every $\tilde{\mathfrak{g}}$-submodule of $N(l\Lambda_{0})$ is an ideal of $N(l\Lambda_{0})$ viewed as a VOA. In 
particular, there exists an induced structure of simple VOA on
$L(l\Lambda_{0})$. 
\end{theo+}  

Any automorphism $\sigma$ of order
$T$ of $\mathfrak{g}$ preserves the form $\langle \cdot,\!\cdot \rangle$ 
and induces a Lie algebra automorphism of $\tilde{\mathfrak{g}}$. It
follows from the associator formula for VOAs that $\sigma$ 
also induces
VOA automorphisms of $N(l\Lambda_{0})$ and $L(l\Lambda_{0})$
respectively, and we denote
these induced automorphisms again by $\sigma$. If $M$ is a restricted
$\hat{\mathfrak{g}}[\sigma]$-module of level $l$, the map
\begin{gather*}
  Y^{\sigma}_{M}:\, \mathfrak{g}(-1)\otimes 1 \longrightarrow (\mbox{End}\,
  M)\big[\big[z^{1/T},z^{-1/T}\big]\big]\\
  Y^{\sigma}_{M}(a(-1)\otimes 1,z)=a(\sigma;z),
\end{gather*}
for $a\in \mathfrak{g}_{(j)}$, $j=0,\ldots,T-1$, has a unique extension
to $N(l\Lambda_{0})$ that makes $(M,Y^{\sigma}_{M})$ a weak
$\sigma$-twisted $N(l\Lambda_{0})$-module. This is a consequence of
the theory of local systems of twisted vertex operators developed in
\cite{Li2}, where the following $\sigma$-twisted counterpart of Theorem
1.1 was obtained:}
{\allowdisplaybreaks
\begin{theo+}
  Let $l\neq -h\check{}$ be a complex number. Then any restricted
  $\hat{\mathfrak{g}}[\sigma]$-module of level $l$ is a weak
  $\sigma$-twisted $N(l\Lambda_{0})$-module. 
\end{theo+}}

\section{Main results}

We use the setting of Section 1 throughout: $\mathfrak{g}$ is a
finite-dimensional simple Lie algebra with the form $\langle \cdot\,
,\cdot \rangle$
normalized as in \S 1.1, $\sigma$ is an automorphism of order $T$ of
$\mathfrak{g}$, $\mathfrak{g}=\coprod_{j\in \mathbf{Z}_{T}}
\mathfrak{g}_{(j)}$ denotes
the $\sigma$-gradation of $\mathfrak{g}$, and $M$ is a restricted
$\tilde{\mathfrak{g}}[\sigma]$-module of level $k\in \mathbf{C}$ (in
particular, $M$ could be a Verma module). Recall from Theorem 1.2 
that $(M, Y_{M}^{\sigma})$ is a weak $\sigma$-twisted
$N(k\Lambda_{0})$-module and let $R$ be a $\sigma$-invariant subspace of
$N(k\Lambda_{0})$ with $\sigma$-decomposition $R=\coprod_{j\in \mathbf{Z}_{T}}
R^{j}$. We also assume that $R$ is invariant under both $\mathfrak{g}(0)$ and 
$L(0)$ ($=\mbox{Res}_{z}zY(\omega,z)$). Define the space 
\begin{equation}\label{2.1}
  \bar{R}_{\sigma}=\text{{\bf C}-span} \, \{r_{n}\mid \, r\in R,\, n\in
  \tfrac{1}{T}\mathbf{Z}\}\subset \mbox{End {\it M}},
\end{equation}
where $Y_{M}^{\sigma}(r,z)=\sum_{n\in \tfrac{1}{T}\mathbf{Z}}r_{n}z^{-n-1}$, 
and let $\eta=\exp(2\pi i/T)$. We denote by
$\sigma$ as well the linear automorphism of $\mbox{(End {\it
    M})}\big[\big[z^{1/T},z^{-1/T}\big]\big]$ determined by $\sigma
f(z^{1/T})=f(\eta^{-1}z^{1/T})$. It follows from
\eqref{1.8} that $\sigma (v_{n})=(\sigma(v))_{n}$ for $v\in N(k\Lambda_{0})$,
$n\in \frac{1}{T}\mathbf{Z}$, so that in particular
$\bar{R}_{\sigma}$ is $\sigma$-stable and thus 
\begin{equation}\label{2.2}
  \bar{R}_{\sigma}=\coprod_{j\in \mathbf{Z}_{T}}\bar{R}_{\sigma}^{j},
\end{equation}
where $\bar{R}_{\sigma}^{j}$ is the $\eta^{j}$-eigenspace of $\sigma$
in $\bar{R}_{\sigma}$. Furthermore, the twisted commutator formula
\eqref{1.9} together with the derivation property \eqref{1.7} yield
\begin{equation}\label{2.3}
  [L(0), Y_{M}^{\sigma}(v,z)]=Y_{M}^{\sigma}(L(0)v,z)+z\frac{d}{dz}Y_{M}^
{\sigma}(v,z)
\end{equation}
for $v\in N(k\Lambda_{0})$, with $L(0)=\mbox{Res}_{z}zY_{M}^{\sigma}
(\omega,z)$ in the left-hand side. Note that we may 
identify $\tilde{\mathfrak{g}}[\sigma]$ with $\hat{\mathfrak{g}}[\sigma]\rtimes
\mathbf{C}L(0)\subset \mbox{End {\it M}}$. Since $R$ is
$L(0)$-invariant, \eqref{2.3} implies that $L(0)$ induces a
$\frac{1}{T}\mbox{{\bf Z}}$-gradation on $\bar{R}_{\sigma}$:
\begin{equation}\label{2.4}
  \bar{R}_{\sigma}=\coprod_{n\in \tfrac{1}{T}\mathbf{Z}}\bar{R}_{\sigma}(n),
\end{equation}
where $\bar{R}_{\sigma}(n)=\{p\in \bar{R}_{\sigma}\mid \,
[L(0),p]=np\}$ for $n\in \frac{1}{T}\mathbf{Z}$. The 
gradations \eqref{2.2} and \eqref{2.4} are compatible in the sense that
\begin{equation}\label{2.5}
  \bar{R}_{\sigma}^{j}=\!\!\coprod_{n\in \tfrac{1}{T}\mathbf{Z},\, nT \bmod
    T\equiv j}\!\!\bar{R}_{\sigma}(n) \text{ for } j\in \mathbf{Z}_{T}\,
  \text{ and } \,
  \bar{R}_{\sigma}(n)=\coprod_{j\in
    \mathbf{Z}_{T}}\bar{R}_{\sigma}(n)^{j} \text{ for } n\in
    \tfrac{1}{T}\mathbf{Z},
\end{equation}
where $\bar{R}_{\sigma}(n)^{j}=\bar{R}_{\sigma}^{j}\cap
\bar{R}_{\sigma}(n)$. Note that \eqref{1.8} implies that
$\bar{R}_{\sigma}(n)$ is in fact $\sigma$-homogeneous, since
$\bar{R}_{\sigma}(n)^{j}=0$ unless $j\equiv nT \bmod T$.
\begin{theorem}
  {\em (i)} Let $R$ and $M$ be as above and assume that $\mathfrak{g}(n)R=0$
  for $n\in \mathbf{Z}_{_{>0}}$. Then
  \begin{equation}\label{2.6}
    [x(m),r_{n}]=(x(0)r)_{m+n}
  \end{equation}
  for all $x\in \mathfrak{g}$, $r\in R$, and $m,n\in \frac{1}{T}\mathbf{Z}$, 
so that $\bar{R}_{\sigma}$ becomes a loop module under the
  adjoint action of $\tilde{\mathfrak{g}}[\sigma]$. Conversely, if $M$ is
  a faithful weak $\sigma$-twisted $N(k\Lambda_{0})$-module and \eqref{2.6}
  holds, then
  \begin{equation}\label{2.7}
    \mathfrak{g}(n)R=0 \quad \forall n\in \mathbf{Z}_{_{>0}}.
  \end{equation}
  {\em (ii)} $\bar{R}_{\sigma}$ is an irreducible loop
  $\tilde{\mathfrak{g}}[\sigma]$-module if $R$ 
  is a nontrivial irreducible $\mathfrak{g}$-module. Moreover, if $M$ is a
  faithful weak $\sigma$-twisted $N(k\Lambda_{0})$-module, then the
  converse is also true.
\end{theorem}  
\begin{proof}
  (i) Suppose that $\mathfrak{g}(n)R=0$ for all positive integers
$n$. It suffices to prove \eqref{2.6} for $\sigma$-homogeneous elements
$x\in \mathfrak{g}_{(j)}$, $r\in R^{k}$, where $j,k\in
\{0,1,\ldots,T-1\}$. Note first that $x(0)r\in R^{j+k}$ and let $m\in
\frac{j}{T}+\mathbf{Z}$, $n\in
\frac{k}{T}+\mathbf{Z}$ be fixed. Then \eqref{1.8} and \eqref{1.9} imply 
that 
\begin{equation}\label{2.8}
  [x(\sigma;z_{1}),Y_{M}^{\sigma}(r,z_{2})]=\sum_{p\in \mathbf{Z},\, q\in 
\frac{j+k}{T}+\mathbf{Z}}(x(0)r)_{q}z_{1}^{-p-\frac{j}{T}-1}z_{2}^
{-q+p+\frac{j}{T}-1}.
\end{equation}
On the other hand
\begin{equation}\label{2.9}
  [x(\sigma;z_{1}),Y_{M}^{\sigma}(r,z_{2})]=\sum_{s\in \frac{j}{T}+\mathbf{Z},
\, t\in \frac{k}{T}+\mathbf{Z}}[x(s),r_{t}]z_{1}^{-s-1}z_{2}^{-t-1}
\end{equation}
by \eqref{1.8}, and then \eqref{2.6} follows by comparing the coefficients of
$z_{1}^{-m-1}z_{2}^{-n-1}$ in the right-hand sides of \eqref{2.8} and 
\eqref{2.9} respectively.

Assume now that $M$ is a faithful weak $\sigma$-twisted
$N(k\Lambda_{0})$-module and that \eqref{2.6} holds, and let $N\in
\mbox{{\bf Z}}_{_{\ge 0}}$ be such that $x(p)r=0$ for all $p\ge N+1$. Then
\eqref{2.6} and \eqref{2.9} yield
$$[x(\sigma;z_{1}),Y_{M}^{\sigma}(r,z_{2})]=z_{1}^{-1}\delta\!\left(\!
\frac{z_{2}}{z_{1}}\!\right)\!\left(\!\frac{z_{2}}{z_{1}}\!\right)
^{\frac{j}{T}}Y_{M}^{\sigma}(x(0)r,z_{2}).$$
Thus
$$\sum_{i=1}^{N}\frac{1}{i!}\!\left(\!\left(\!\frac{\partial}{\partial
  z_{2}}\!\right)^{i}z_{1}^{-1}\delta\!\left(\!\frac{z_{2}}{z_{1}}\!\right)\!
\left(\!\frac{z_{2}}{z_{1}}\!\right)^{\frac{j}{T}}\!\right)Y_{M}^{\sigma}
(x(i)r,z_{2})=0$$
by \eqref{1.9}. Using \cite[Lemma 2.3]{Li2} one gets that
$Y_{M}^{\sigma}(x(i)r,z_{2})=0$ for $1\le i\le N$. Therefore $x(i)r=0$
for $i\in \{1,\ldots,N\}$ as well (since $M$ is faithful), which
proves \eqref{2.7}.
\newline
(ii) We first prove the second statement. Suppose that $\bar{R}_{\sigma}$ is 
irreducible. If $R$ were a
1-dimensional trivial $\mathfrak{g}$-module, then $\mbox{{\bf C}}r_{n}$
would be 
a nonzero $\tilde{\mathfrak{g}}[\sigma]$-submodule of $\bar{R}_{\sigma}$
for each nonzero $r\in R$ and $n\in \frac{1}{T}\mbox{{\bf Z}}$ satisfying 
$r_{n}\neq 0$ (such elements exist since $M$ is faithful), hence
a contradiction. Suppose now that $S\subset R$ is a nonzero proper
$\mathfrak{g}$-submodule of $R$. Then there exists $0\neq r\in
R_{d}\setminus S_{d}$ for some $d\in \mbox{{\bf Z}}$, where
$R=\coprod_{n\in \mathbf{Z}}R_{n}$ and $S=\coprod_{n\in \mathbf{Z}}S_{n}$ are 
the $L(0)$-gradations of $R$ and $S$
respectively. Note that $\mbox{($\overline{S_{d}}$)$_{\sigma}$}\subset
\mbox{($\overline{R_{d}}$)$_{\sigma}$}$ are
$\tilde{\mathfrak{g}}[\sigma]$-submodules of $\bar{R}_{\sigma}$ by \eqref{2.6},
and thus
\begin{equation}\label{2.10}
  \mbox{($\overline{S_{d}}$)$_{\sigma}$}=
\mbox{($\overline{R_{d}}$)$_{\sigma}$}=\bar{R}_{\sigma}
\end{equation}
since $\bar{R}_{\sigma}$ is irreducible. It follows from \eqref{2.3}, 
\eqref{2.4} and \eqref{2.10} that there exist 
$m\in \frac{1}{T}\mbox{{\bf Z}}$ and
$s\in S_{d}$ such that $r_{m}\neq 0$ and $r_{m}=s_{m}$. Then \eqref{2.6} 
implies that $x(0)t=0$ for all $x\in \mathfrak{g}$, where $t:=r-s$. 
Consequently, $\bar{T}_{\sigma}:=\mbox{{\bf
    C}-span}\,\{t_{l}\mid l\in \frac{1}{T}\mbox{{\bf Z}}\}$ is a
(trivial) $\tilde{\mathfrak{g}}[\sigma]$-submodule of $\bar{R}_{\sigma}$
by \eqref{2.6}. Clearly, $r_{m}\notin \bar{T}_{\sigma}$, so that
$\bar{T}_{\sigma}\neq \bar{R}_{\sigma}$ and thus
$\bar{T}_{\sigma}=0$. Hence $t=0$ (again by the faithfulness of
$M$), and one gets that $r=s\in S_{d}$, which is a contradiction. $R$
must therefore be an irreducible $\mathfrak{g}$-module.

Conversely, suppose that $R$ is a nontrivial irreducible
$\mathfrak{g}$-module. Then $R$ is $L(0)$-homogeneous of some conformal
weight $d\in \mbox{{\bf Z}}$. Recall \eqref{2.5} and notice that
\begin{equation}\label{2.11}
  \bar{R}_{\sigma}(n)=\{r_{d-n-1}\mid r\in R\},\quad n\in
  \tfrac{1}{T}\mbox{{\bf Z}}.
\end{equation}
Let $Z$ be a nonzero $\tilde{\mathfrak{g}}[\sigma]$-submodule of
$\bar{R}_{\sigma}$. Then $Z=\coprod_{n\in \frac{1}{T}\mathbf{Z}}Z(n)$ 
($Z$ being $L(0)$-stable), where $Z(n)=\{z\in Z\mid
[L(0),z]=nz\}\subset \bar{R}_{\sigma}(n)$. By \eqref{2.11} there exist $n\in
  \frac{1}{T}\mbox{{\bf Z}}$ and $r\in R$ such that $0\neq
  r_{d-n-1}\in Z(n)$, and it follows from \eqref{2.6} that
\begin{equation}\label{2.12}
  (x(0)r)_{d-m-1}=[x(n-m),r_{d-n-1}]\in Z(m)
\end{equation}
for every $x\in \mathfrak{g}$ and $m\in \frac{1}{T}\mbox{{\bf
      Z}}$. Since $R$ is irreducible, one has
  $R=U(\mathfrak{g})\cdot r$, and then by iterating \eqref{2.12} one gets
  $\bar{R}_{\sigma}(m)\subset Z(m)$ for all $m\in \frac{1}{T}\mbox{{\bf
      Z}}$. Therefore $Z=\bar{R}_{\sigma}$ and the proof is
  complete.
\end{proof}  

\begin{remark}
It follows from the twisted associator formula \eqref{1.10} 
and the irreducibility
of $L(l\Lambda_{0})$ that level $l$ standard
$\hat{\mathfrak{g}}[\sigma]$-modules are in fact faithful
$\sigma$-twisted $L(l\Lambda_{0})$-modules. It would be interesting to see 
whether the Verma $\tilde{\mathfrak{g}}[\sigma]$-module $M(\Lambda)$ is 
itself a faithful $\sigma$-twisted $N(k\Lambda_{0})$-module, where 
$\Lambda\in P_{+}$ is such that $\Lambda(c)=k\in \mathbf{Z}_{+}$. This 
question is obviously related to the classification of submodules of 
$N(k\Lambda_{0})$, which in turn lies close to the problem of classifying 
annihilating ideals of standard (or even admissible) $\tilde{\mathfrak{g}}
[\sigma]$-modules (cf.\@~\cite{FM}).
\end{remark}

{\em In the remainder of this section, $k$ is assumed to be a fixed
  positive integer.} The next result will be used in the proofs of
Theorems 2.9 and 2.13.
\begin{proposition}[\cite{Li2}]
  Let $V$ be a \mbox{{\em VOA}}, let $\sigma$ be an
  automorphism of order $S$ of $V$ and $h\in V$ such that 
  $L(n)h=\delta_{n,0}h$, $\sigma(h)=h$, $[h_{m},h_{n}]=0$  
  for $m,n\in \mathbf{Z}_{_{\ge 0}}$. 
  Assume that $h(0)$ acts semisimply on $V$ and that
  $\mbox{{\em Spec}$(h(0))\subset \frac{1}{T}${\em {\bf Z}}}$ for some $T\in
  \mbox{{\em {\bf Z}}}_{_{>0}}$, so that $\sigma_{h}:=\exp(2\pi ih(0))$ is
  an automorphism of $V$ satisfying
  $\sigma_{h}^{T}=\mbox{{\em id}}_{_{V}}$. Let $(M, Y_{M}^{\sigma})$ be a
  $\sigma$-twisted $V$-module and define
  \begin{equation}\label{2.13}
    \Delta(h,z)=z^{h(0)}\exp\left(\sum_{n=1}^{\infty}(-1)^{n-1}\frac{h(n)}{n}
z^{-n}\right),\quad \overline{Y}_{M}^{\sigma}(a,z)=Y_{M}^{\sigma}
(\Delta(h,z)a,z)
  \end{equation}
  for $a\in V$. Then $(M,\overline{Y}_{M}^{\sigma}(\cdot,z))$ is a weak
  $(\sigma \sigma_{h})$-twisted $V$-module.
\end{proposition}

\begin{remark}
  As pointed out in \cite{Li2}, Proposition 2.3 gives an
  isomorphism between 
  $\hat{\mathfrak{g}}[\sigma]$ and $\hat{\mathfrak{g}}[\sigma \sigma_{h}]$. 
  Furthermore, $(M,\overline{Y}_{M}^{\sigma})$ is
  irreducible if $(M,Y_{M}^{\sigma})$ is
  irreducible ($\Delta(h,z)$ being invertible), and any
  $\sigma_{h}$-twisted $V$-module can be constructed from a
  $V$-module.
\end{remark}

The following well-known result will be needed for the proof of Theorem 2.13:
\begin{proposition}[\cite{K}]
  Let $\mathfrak{g}$ be a simple finite-dimensional Lie algebra,
  let $\mathfrak{t}$ be a \mbox{{\em CSA}} of $\mathfrak{g}$ and let
  $\Pi=\{\alpha_{1},\ldots,\alpha_{n}\}$ be a set of simple roots. Let
  $\sigma\in \mbox{{\em Aut}}(\mathfrak{g})$ be such that
  $\sigma^{T}=\mbox{{\em id}}_{\mathfrak{g}}$. Then $\sigma$ is conjugate
  to an automorphism of $\mathfrak{g}$ of the form
  \begin{equation}\label{2.14}
    \mu \exp\!\big(\mbox{{\em ad}}\!\left(\!\tfrac{2\pi
      i}{T}h\!\right)\!\big),\quad h\in \mathfrak{t}_{[0]},
  \end{equation}
  where $\mu$ is a diagram automorphism preserving $\mathfrak{t}$ and
  $\Pi$, $\mathfrak{t}_{[0]}$ is the fixed-point set of $\mu$ in
  $\mathfrak{t}$, and $\alpha_{i}(h)\in \mbox{{\em {\bf Z}}}$ for $i\in
  \{1,\ldots,n\}$.
\end{proposition}
As in \S 1.1, we fix a $CSA$ of $\mathfrak{g}$ and denote it by
$\mathfrak{t}$. Let $\Phi$ be the root system of $\mathfrak{g}$ and
let $\Pi=\{\alpha_{1},\ldots,\alpha_{n}\}\subset \mathfrak{t}^{*}$ be
a basis of $\Phi$ enumerated as in \cite[Table {\em Fin}]{K}. Choose root 
vectors $x_{\alpha}\in
\mathfrak{g}_{\alpha}$ such that $h_{\alpha}:=[x_{\alpha},x_{-\alpha}]$
satisfies $\alpha(h_{\alpha})=2$ for $\alpha\in \Phi$, and let $\theta$
be the highest root. Notice that with the above normalizations one
always has $\langle x_{\theta}, x_{-\theta}\rangle =1$. We shall use the 
following
result from the untwisted representation theory (cf.\@~\cite{Li1,MP2}):
\begin{proposition}
  Let $M$ be an integrable $\tilde{\mathfrak{g}}$-module of level
  $k$. Then $x_{\alpha}(z)^{tk+1}=0$ acting on $M$, where $t=1$ if
  $\alpha$ is a long root, $t=2$ if $\alpha$ is a short root and
  $\mathfrak{g}$ is not of type $G_{2}$, and $t=3$ if $\alpha$ is a
  short root and $\mathfrak{g}$ is of type $G_{2}$.
\end{proposition}
Let now $\mu$ be a diagram automorphism of $\mathfrak{g}$ induced by
an automorphism $\bar{\mu}$ of the Dynkin diagram of $\mathfrak{g}$
with respect to $(\mathfrak{t},\Pi)$. It is well-known that the
subgroup of $\mbox{Aut}(\mathfrak{g})$ generated by all such diagram
automorphisms is isomorphic to the symmetric group $\mathcal{S}_{m}$, 
where $m=1$ for $B_{n}$, $C_{n}$, $E_{7}$, $E_{8}$, $F_{4}$ and
$G_{2}$, $m=2$ for $A_{n}$ ($n\ge 2$), $D_{n}$ ($n\ge 3$, $n\neq 4$) 
and $E_{6}$, and $m=3$ for $D_{4}$. Let $r$ ($=1,2,$ or $3$)
be the order of $\mu$, so that in particular one has the
$\mu$-decompositions $\mathfrak{g}=\coprod_{j\in 
  \mbox{{\bf Z}}_{r}}\mathfrak{g}_{[j]}$ and (since $\mathfrak{t}$
is obviously $\mu$-stable) $\mathfrak{t}=\coprod_{j\in
  \mbox{{\bf Z}}_{r}}\mathfrak{t}_{[j]}$. We shall concentrate on the
five cases when $r\ge 2$, namely when $\mathfrak{g}$ is of type $A_{2l}$,
$A_{2l-1}$, $D_{l+1}$, $E_{6}$ or $D_{4}$, and $r=2,2,2,2 \text{ or }
3$ respectively (for $r=3$ there are two
  equivalent automorphisms of this type and we just choose one of
  them). Recall that $\mathfrak{g}_{[0]}$ is a simple
subalgebra of $\mathfrak{g}$ such that $\mathfrak{t}_{[0]}=\mathfrak{g}_{[0]}
\cap \mathfrak{t}$ is a $CSA$ of
$\mathfrak{g}_{[0]}$. More specifically, $\mathfrak{g}_{[0]}$ is of
type $B_{l}$, $A_{1}$, $C_{l}$, $B_{l}$, $F_{4}$, $G_{2}$ whenever
$\mathfrak{g}$ is of
type $A_{2l}$ ($l\ge 2$), $A_{2}$, $A_{2l-1}$, $D_{l+1}$, $E_{6}$,
$D_{4}$ respectively. Recall also the
elements $E_{i},F_{i},H_{i}\in \mathfrak{g}$ and the simple roots
$\beta_{i}\in \mathfrak{t}_{[0]}^{*}$, $i=0,1,\ldots l$, defined in \S
1.1. These are given explicitly in \cite{K}, 
where an additional element  
$\theta^{0}\in \Phi$ is introduced as follows (in case 1 we switched
the indexes $0$ and $l$ as compared with \cite{K}):
{\allowdisplaybreaks 
\begin{align*}
  \mbox{{\em Case 1:}}\quad & \mathfrak{g}=A_{2l}, \, r=2; \,
   \bar{\mu}(\alpha_{i})=\alpha_{2l-i+1} \, (1\le i\le 2l), \,
  \mu(h_{\alpha})=h_{\bar{\mu}(\alpha)}, \\
  & \mu(x_{\alpha})=(-1)^{1+ht(\alpha)}x_{\bar{\mu}(\alpha)};\,
  \theta=\sum_{i=1}^{2l}\alpha_{i}=\theta^{0}, \,
  -\beta_{0}=\theta^{0}|_{\mathfrak{t}_{[0]}}, \\
  & \beta_{i}=\alpha_{i}|_{\mathfrak{t}_{[0]}}=\alpha_{2l-i+1}|_{\mathfrak{t}
_{[0]}} \, (1\le i\le l); \, E_{0}=x_{-\theta}, \\
  & E_{i}=x_{\alpha_{i}}+x_{\alpha_{2l-i+1}} 
  (1\le i\le l-1), \, E_{l}=\sqrt{2}(x_{\alpha_{l}}+x_{\alpha_{l+1}});\, 
  F_{0}=x_{\theta}, \\
  & F_{i}=x_{-\alpha_{i}}+x_{-\alpha_{2l-i+1}}
  \, (1\le i\le l-1), \,
  F_{l}=\sqrt{2}(x_{-\alpha_{l}}+x_{-\alpha_{l+1}}); \\
  & H_{i}=[E_{i},F_{i}] \, (0\le i\le l). \\ 
\mbox{{\em Case 2:}} \quad & \mathfrak{g}=A_{2l-1}, \, r=2; \,
  \bar{\mu}(\alpha_{i})=\alpha_{2l-i} \: (1\le i\le 2l-1), \,
  \mu(h_{\alpha})=h_{\bar{\mu}(\alpha)}, \\
  & \mu(x_{\alpha})=x_{\bar{\mu}(\alpha)};\, \theta^{0}=\theta-\alpha_{2l-1}, \,
  -\beta_{0}=\frac{1}{2}\big(\theta^{0}+\bar{\mu}(\theta^{0})\big)\big|_
{\mathfrak{t}_{[0]}}, \\
  & \beta_{i}=\alpha_{i}|_{\mathfrak{t}_{[0]}}=\alpha_{2l-i}|_{\mathfrak{t}
_{[0]}} \: (1\le i\le l); \, E_{0}=x_{-\theta^{0}}-x_{-\bar{\mu}(\theta^{0})}, \\
  & E_{i}=x_{\alpha_{i}}+x_{\alpha_{2l-i}} \:
  (1\le i\le l-1), \, E_{l}=x_{\alpha_{l}};\,
  F_{0}=x_{\theta^{0}}-x_{\bar{\mu}(\theta^{0})}, \\
  & F_{i}=x_{-\alpha_{i}}+x_{-\alpha_{2l-i}}
  \: (1\le i\le l-1), \,
  F_{l}=x_{-\alpha_{l}}; \,
  H_{i}=[E_{i},F_{i}] \: (0\le i\le l). \\ 
\mbox{{\em Case 3:}} \quad & \mathfrak{g}=D_{l+1}, \: r=2; \,
  \bar{\mu}(\alpha_{i})=\alpha_{i} \:(1\le i\le l-1), \: \bar{\mu}
(\alpha_{l})=\alpha_{l+1}, \:
  \bar{\mu}(\alpha_{l+1})=\alpha_{l}, \\
  & \mu(h_{\alpha})=
h_{\bar{\mu}(\alpha)},\, \mu(x_{\alpha})=x_{\bar{\mu}(\alpha)}; \: 
  \theta=\alpha_{1}+2\sum_{i=2}^{l-1}\alpha_{i}+\alpha_{l}+\alpha_{l+1}, \\ 
  & \theta^{0}=\frac{1}{2}(\theta+\alpha_{1}+\alpha_{l}-\alpha_{l+1}), \, -\beta_{0}=\frac{1}{2}\big(\theta^{0}+\bar{\mu}(\theta^{0})\big)\big|_
{\mathfrak{t}_{[0]}}, \\
  & \beta_{i}=\alpha_{i}|_{\mathfrak{t}_{[0]}} \: 
(1\le i\le l-1), \:
  \beta_{l}=\alpha_{l}|_{\mathfrak{t}_{[0]}}=\alpha_{l+1}|_{\mathfrak{t}
_{[0]}}; \\
  & E_{0}=x_{-\theta^{0}}-x_{-\bar{\mu}(\theta^{0})}, \:
  E_{i}=x_{\alpha_{i}} \: (1\le i\le l-1), \:
  E_{l}=x_{\alpha_{l}}+x_{\alpha_{l+1}}; \\
  & F_{0}=x_{\theta^{0}}-x_{\bar{\mu}(\theta^{0})}, \,
  F_{i}=x_{-\alpha_{i}} \: (1\le i\le l-1), \:
  F_{l}=x_{-\alpha_{l}}+x_{-\alpha_{l+1}}; \\
  & H_{i}=[E_{i},F_{i}] \: (0\le i\le l). \\ 
\mbox{{\em Case 4:}} \quad & \mathfrak{g}=E_{6}, \: r=2; \,
  \bar{\mu}(\alpha_{1})=\alpha_{5}, \, \bar{\mu}(\alpha_{2})=\alpha_{4}, \, 
\bar{\mu}(\alpha_{3})=\alpha_{3}, \,
\bar{\mu}(\alpha_{6})=\alpha_{6}, \\ 
  & \mu(h_{\alpha})=h_{\bar{\mu}(\alpha)},\, 
  \mu(x_{\alpha})=x_{\bar{\mu}(\alpha)}; \:
   \theta=\alpha_{1}+2\alpha_{2}+3\alpha_{3}+2\alpha_{4}+\alpha_{5}+
2\alpha_{6}, \\
  & \theta^{0}=\theta-\alpha_{3}-\alpha_{4}-\alpha_{6}, \,
  -\beta_{0}=\frac{1}{2}\big(\theta^{0}+\bar{\mu}(\theta^{0})\big)\big|
_{\mathfrak{t}_{[0]}}, \, 
  \beta_{1}=\alpha_{1}|_{\mathfrak{t}_{[0]}}=\alpha_{5}|_{\mathfrak{t}_{[0]}},
  \\
  & \beta_{2}=\alpha_{2}|_{\mathfrak{t}_{[0]}}=\alpha_{4}|_{\mathfrak{t}_{[0]}},
  \, 
  \beta_{3}=\alpha_{3}|_{\mathfrak{t}_{[0]}}, \,
  \beta_{4}=\alpha_{6}|_{\mathfrak{t}_{[0]}}; \, 
  E_{0}=x_{-\theta^{0}}-x_{-\bar{\mu}(\theta^{0})}, \\
  & E_{1}=x_{\alpha_{1}}+x_{\alpha_{5}}, \,
  E_{2}=x_{\alpha_{2}}+x_{\alpha_{4}}, \, 
  E_{3}=x_{\alpha_{3}}, \\
  & E_{4}=x_{\alpha_{6}}; \, 
  F_{0}=x_{\theta^{0}}-x_{\bar{\mu}(\theta^{0})}, \,
  F_{1}=x_{-\alpha_{1}}+x_{-\alpha_{5}}, \,
  F_{2}=x_{-\alpha_{2}}+x_{-\alpha_{4}}, \\ 
  & F_{3}=x_{-\alpha_{3}}, \,
  F_{4}=x_{-\alpha_{6}}; \, 
  H_{i}=[E_{i},F_{i}] \: (0\le i\le 4). \\ 
\mbox{{\em Case 5:}} \quad & \mathfrak{g}=D_{4}, \, r=3, \,
  \varepsilon=\exp(2\pi i/3); \,
  \bar{\mu}(\alpha_{1})=\alpha_{4}, \, \bar{\mu}(\alpha_{2})=\alpha_{2}, \, 
\bar{\mu}(\alpha_{3})=\alpha_{1}, \\
  & \bar{\mu}(\alpha_{4})=\alpha_{3}, \, \mu(h_{\alpha})=h_{\bar{\mu}(\alpha)},\,
  \mu(x_{\alpha})=x_{\bar{\mu}(\alpha)}; \:
  \theta=\alpha_{1}+2\alpha_{2}+\alpha_{3}+\alpha_{4}, \\
  & \theta^{0}=\theta-\alpha_{2}-\alpha_{4}, \,
  -\beta_{0}=\frac{1}{3}\big(\theta^{0}+\bar{\mu}(\theta^{0})+\bar{\mu}^{2}
(\theta^{0})\big)\big|_{\mathfrak{t}_{[0]}}, \\
  & \beta_{1}=\alpha_{1}|_{\mathfrak{t}_{[0]}}=\alpha_{3}|_{\mathfrak{t}_{[0]}}
=\alpha_{4}|_{\mathfrak{t}_{[0]}}, \, \beta_{2}=\alpha_{2}|_{\mathfrak{t}
_{[0]}}; \\
  & E_{0}=x_{-\theta^{0}}+\varepsilon^{2}
  x_{-\bar{\mu}(\theta^{0})}+\varepsilon
  x_{-\bar{\mu}^{2}(\theta^{0})}, \,
  E_{1}=x_{\alpha_{1}}+x_{\alpha_{3}}+x_{\alpha_{4}},\,
  E_{2}=x_{\alpha_{2}}; \\
  & F_{0}=x_{\theta^{0}}+\varepsilon
  x_{\bar{\mu}(\theta^{0})}+\varepsilon^{2}x_{\bar{\mu}^{2}(\theta^{0})}, \,
  F_{1}=x_{-\alpha_{1}}+x_{-\alpha_{3}}+x_{-\alpha_{4}}, \,
  F_{2}=x_{-\alpha_{2}};\\
  & H_{i}=[E_{i},F_{i}] \: (0\le i\le 2).
\end{align*}
\indent
Note that $\mathfrak{g}$ is simply-laced and that
$x_{\theta}$ is $\mu$-homogeneous in all these cases.} 
 The canonical generators
$\{e_{i},\, f_{i},\, h_{i}\mid 0\le i\le l\}$ of $\hat{\mathfrak{g}}[\mu]$
 are as in \eqref{1.3}, with $\nu=\mu$ and
$(s_{0},s_{1},\ldots,s_{l})=(1,0,\ldots,0)$. Set now 
\begin{equation}\label{2.15}
  R=U(\mathfrak{g})x_{\theta}(-1)^{k+1}\mbox{{\bf 1}}\subset
  N(k\Lambda_{0}).
\end{equation}
Clearly, $R$ is isomorphic to the simple 
$\mathfrak{g}$-module with highest weight $(k+1)\theta$. 
Notice also that $R$ is invariant under all
automorphisms of the form \eqref{2.14} and recall
from \eqref{1.4} the triangular decomposition of $\hat{\mathfrak{g}}$. The 
following result may be found in \cite[Corollary 5.4 \& Lemma 5.5]{MP2}:
\begin{proposition}
  The maximal submodule $N^{1}(k\Lambda_{0})$ of $N(k\Lambda_{0})$ is
  generated by the singular vector $x_{\theta}(-1)^{k+1}\mbox{{\em {\bf
      1}}}$. Hence 
$$N^{1}(k\Lambda_{0})=U(\hat{\mathfrak{g}})x_{\theta}(-1)^{k+1}\mbox{{\em {\bf
      1}}}=U(\mathfrak{n}_{-})x_{\theta}(-1)^{k+1}\mbox{{\em {\bf 1}}}.$$
\end{proposition}
  Let $M$
  be a restricted $\tilde{\mathfrak{g}}[\mu]$-module of level $k$ and
  define $\bar{R}_{\mu}$ as in \eqref{2.1}. Proposition 2.7 implies that
  $\mathfrak{g}(n)R=0$ for all integers $n\ge 1$, and then Theorem 2.1
  yields
\begin{corollary}
  $\bar{R}_{\mu}$ is an irreducible loop
  $\tilde{\mathfrak{g}}[\mu]$-module.\hfill $\Box$
\end{corollary}
Recall from Theorem 1.2 that any restricted (in particular, any
integrable) $\tilde{\mathfrak{g}}[\mu]$-module of level $k$ is a weak
$\mu$-twisted $N(k\Lambda_{0})$-module. We can now prove the following 
\begin{theorem}
  If $L(\Lambda)$ is a standard $\tilde{\mathfrak{g}}[\mu]$-module of
  level $k$ then $\bar{R}_{\mu}L(\Lambda)=0$.
\end{theorem}
\begin{proof}
  Let $M=L(\Lambda)$. By the twisted associator formula \eqref{1.10} and
  induction, it suffices to prove that
  $Y_{M}^{\mu}(x_{\theta}(-1)^{k+1}\mbox{{\bf 1}},z)=0$. Recall that
  $\langle x_{\theta},x_{-\theta} \rangle =1$ and notice that $x_{\theta}\in
  \mathfrak{g}_{[1]}$ in case 1, while $x_{\theta}\in
  \mathfrak{g}_{[0]}$ in cases 2-5. Let $\mathfrak{a}=\mbox{{\bf
      C}-span}\{x_{\theta},x_{-\theta},h_{\theta}\} \, (\cong
  \mathfrak{sl}(2,\mbox{{\bf C}}))$. Then we can embed
  $\widehat{\mathfrak{sl}(2,\mbox{{\bf C}})}\cong \hat{\mathfrak{a}}$ into
  $\hat{\mathfrak{g}}[\mu]$ with canonical central element $c$ in cases 2-5. 
Since $M$ is a
  level $k$ integrable $\hat{\mathfrak{g}}[\mu]$-module, it is {\it \`{a}
    fortiori} an integrable $\hat{\mathfrak{a}}$-module of level $k$. As
  such, $(M,Y_{M}^{\mu})$ becomes by Theorem 1.1 an (untwisted) weak
  module for the VOA $N(k\Lambda_{0};\hat{\mathfrak{a}})\subset
  N(k\Lambda_{0})$ via the restricted vertex map. It then follows from
  Proposition 2.6 that $Y_{M}^{\mu}(x_{\theta}(-1)\mbox{{\bf
      1}},z)^{k+1}=0$. Since $[Y_{M}^{\mu}(x_{\theta}(-1)\mbox{{\bf
      1}},z_{1}),Y_{M}^{\mu}(x_{\theta}(-1)\mbox{{\bf
      1}},z_{2})]=0$, one gets from \eqref{1.10} that 
$$Y_{M}^{\mu}(x_{\theta}
(-1)^{k+1}\mbox{{\bf
      1}},z)=Y_{M}^{\mu}(x_{\theta}(-1)\mbox{{\bf
      1}},z)^{k+1}=0$$ 
in cases 2-5, as needed.

  Note that in case 1 we can embed
  $\widehat{\mathfrak{sl}(2,\mbox{{\bf 
      C}})}[\mu]\cong \hat{\mathfrak{a}}[\mu]$ into
  $\hat{\mathfrak{g}}[\mu]$ with canonical central element $c$ and
  that the restriction of 
  $\mu$ to $N(k\Lambda_{0};\hat{\mathfrak{a}})$ coincides with the inner 
  automorphism $\exp(\pi ih_{\theta}(0)/2)$. As above, $M$ is an
  integrable $\hat{\mathfrak{a}}[\mu]$-module of level $k$ and thus
  $(M,Y_{M}^{\mu})$ becomes by Theorem 1.2 a weak $\mu$-twisted
  $N(k\Lambda_{0};\hat{\mathfrak{a}})$-module via the restricted vertex
  map. Using Proposition 2.3 with $h=-h_{\theta}/4$ and
  $\overline{Y}_{M}^{\mu}(\cdot\,,z)=Y_{M}^{\mu}(\Delta(h,z)\cdot\,,z)$, one 
gets that $(M,\overline{Y}_{M}^{\mu})$ is an untwisted weak
  $N(k\Lambda_{0};\hat{\mathfrak{a}})$-module. Consequently, $M$
  acquires a structure of level $k$ integrable
  $\hat{\mathfrak{a}}$-module via the vertex map
  $\overline{Y}_{M}^{\mu}(\cdot\,,z)$. Moreover, it follows from 
$h_{\theta}(n)x_{\theta}(-1)^{p}\mbox{{\bf
      1}}=0$ for $n,p\in
  \mbox{{\bf Z}}_{_{>0}}$ that $\Delta(h,z)x_{\theta}(-1)^{p}\mbox{{\bf
      1}}=z^{-p/2}x_{\theta}(-1)^{p}\mbox{{\bf
      1}}$, which combined with \eqref{1.10} and \eqref{2.13} gives 
  \begin{align}
    \overline{Y}_{M}^{\mu}(x_{\theta}(-1)^{p}\mbox{{\bf 1}},z)=z^{-p/2}
{Y}_{M}^{\mu}(x_{\theta}&(-1)^{p}\mbox{{\bf 1}},z)\label{2.16}\\
&=z^{-p/2}
{Y}_{M}^{\mu}(x_{\theta}(-1)\mbox{{\bf
    1}},z)^{p}=\overline{Y}_{M}^{\mu}(x_{\theta}(-1)\mbox{{\bf
    1}},z)^{p}.\nonumber
  \end{align}
Proposition 2.6 again implies that $\overline{Y}_{M}^{\mu}
(x_{\theta}(-1)^{k+1}\mbox{{\bf
    1}},z)=\overline{Y}_{M}^{\mu}(x_{\theta}(-1)\mbox{{\bf
    1}},z)^{k+1}=0$, and then \eqref{2.16} yields ${Y}_{M}^{\mu}
(x_{\theta}(-1)^{k+1}\mbox{{\bf
    1}},z)=0$. The proof is thereby complete.
\end{proof}

\smallskip

Since $\mathfrak{g}$ is simply-laced in all the cases listed above, one gets 
from \eqref{2.15} and the action of the Weyl group of $\mathfrak{g}$ that 
\begin{equation}\label{2.17}
  \mbox{$x_{\alpha}(-1)^{k+1}\mbox{{\bf 1}}\in R$ \, for 
  $\alpha\in \Phi$}.
\end{equation}

Dealing with case 1 requires the following
\begin{lemma}
  Let $\mathfrak{l}$ be a 3-dimensional Heisenberg algebra with basis
  $\{x,y,z\}$ such that $[x,y]=z$, $[x,z]=[y,z]=0$. For every
  $m\in \mathbf{Z}_{_{\ge 0}}$, the following holds in $U(\mathfrak{l}):$
  \begin{equation}\label{2.18}
    (x+y)^{m}\in \sum_{
    a+b+2c=m \atop
    a\le b}\mathbf{C}x^{a}y^{b}z^{c}+\sum_{
    a+b+2c=m\atop
    b\le a}\mathbf{C}y^{b}x^{a}z^{c}.
  \end{equation}
\end{lemma}
\begin{proof}
  Let $\omega$ be the associative algebra automorphism of
  $U(\mathfrak{l})$ induced by the Lie algebra involution of $\mathfrak{l}$
  defined by $\omega(x)=y$, $\omega(y)=x$, $\omega(z)=-z$. Let
  $a,b,c\in \mathbf{Z}_{_{\ge 0}}$ be such that $a\le b$ and $a+b+2c=n$, where
  $n\in \mathbf{Z}_{_{\ge 0}}$ is arbitrarily fixed. An easy induction gives
  \begin{equation}\label{2.19}
    xy^{p}=y^{p}x+py^{p-1}z,\quad yx^{p}=x^{p}y-px^{p-1}z
  \end{equation}
  for $p\ge 1$, so that if $a<b$ it follows that 
  \begin{equation}\label{2.20}
    x^{a}y^{b}z^{c}x\in \sum_{
    \alpha + \beta +2\gamma =n+1\atop
    \alpha \le \beta}\mbox{{\bf C}}x^{\alpha}y^{\beta}z^{\gamma}.
  \end{equation}
  Suppose now that $a=b$. Using \eqref{2.19} and induction one
  gets that $x^{p}y^{p}\in \sum_{\alpha +\gamma =n}\mbox{{\bf
        C}}y^{\alpha}x^{\alpha}z^{\gamma}$ 
  whenever $p\ge 1$, and therefore
  \begin{equation}\label{2.21}
    x^{a}y^{a}z^{c}x\in \sum_{
    \alpha + \beta +2\gamma =n+1\atop 
    \beta \le \alpha}\mbox{{\bf C}}y^{\beta}x^{\alpha}z^{\gamma}.
  \end{equation}
  It follows from \eqref{2.20} and \eqref{2.21} that in both cases
  \begin{equation}\label{2.22}
    x^{a}y^{b}z^{c}x\in \sum_{
    \alpha + \beta +2\gamma =n+1\atop
    \alpha \le \beta}\mbox{{\bf C}}x^{\alpha}y^{\beta}z^{\gamma}+\sum_{
    \alpha + \beta +2\gamma =n+1\atop
    \beta \le \alpha}\mbox{{\bf C}}y^{\beta}x^{\alpha}z^{\gamma}.
  \end{equation}
Consequently
\begin{equation}\label{2.23}
    y^{a}x^{b}z^{c}y\in \sum_{
    \alpha + \beta +2\gamma =n+1\atop
    \alpha \le \beta}\mbox{{\bf C}}y^{\alpha}x^{\beta}z^{\gamma}+\sum_{
    \alpha + \beta +2\gamma =n+1\atop
    \beta \le \alpha}\mbox{{\bf C}}x^{\beta}y^{\alpha}z^{\gamma}
  \end{equation}
by using \eqref{2.22} and the action of $\omega$. Then \eqref{2.18} follows by
induction on $m$ from
$(x+y)^{m}=(x+y)^{m-1}(x+y)$ together with \eqref{2.22} and \eqref{2.23}.
\end{proof}
\begin{proposition}
  Let $\mathfrak{g}$ be as above. Then $F_{i}(-1)^{n}\mbox{{\em {\bf 1}}}\in
  N^{1}(k\Lambda_{0})$ for all $i\in \{0,1,\ldots,l\}$ if $n$ is
  sufficiently large.
\end{proposition}    
\begin{proof}
  Using the explicit form of the root systems $\Phi$, it is readily
  checked that in cases 2-5 the $F_{i}$'s are sums of at most three 
  commuting root vectors of $\mathfrak{g}$. Therefore
  $$F_{i}(-1)^{tk+1}\mbox{{\bf 1}}\in
  N^{1}(k\Lambda_{0}),\quad i\in \{0,1,\ldots,l\},$$
for some $t\in \{1,2,3\}$ by \eqref{2.17}. The same argument implies 
that in case 1
$$F_{i}(-1)^{2k+1}\mbox{{\bf 1}}\in
  N^{1}(k\Lambda_{0}),\quad i\in \{1,\ldots,l-1\},$$
while obviously
$F_{0}(-1)^{k+1}\mbox{{\bf 1}}=x_{\theta}(-1)^{k+1}\mbox{{\bf 1}}\in
N^{1}(k\Lambda_{0})$. Applying now Lemma 2.10 with $x=x_{-\alpha_{l}}(-1)$,
$y=x_{-\alpha_{l+1}}(-1)$, and
$z=[x_{-\alpha_{l}}(-1),x_{-\alpha_{l+1}}(-1)]$ one gets 
$F_{l}(-1)^{4k+1}\mbox{{\bf 1}}\in N^{1}(k\Lambda_{0})$, 
so that the proposition is true in all cases.
\end{proof}

\smallskip

We are now ready to prove a more complete version of Theorem 2.9:
\begin{theorem}
  Let $M(\Lambda)$ denote the Verma $\hat{\mathfrak{g}}[\mu]$-module with
  highest weight $\Lambda\in P_{+}$ such that $\Lambda(c)=k$. Then
  $\bar{R}_{\mu}M(\Lambda)=M^{1}(\Lambda)$, 
where $M^{1}(\Lambda)$ is the maximal submodule of $M(\Lambda)$.
\end{theorem}
\begin{proof}
  Let $v_{\Lambda}$ be a highest weight vector of $M:=M(\Lambda)$ and
  set $V=\bar{R}_{\mu}M(\Lambda)$. By Theorem 2.9 it suffices to prove
  that $M^{1}(\Lambda)\subset V$, which in turn reduces to showing
  that
  $$\{f_{i}^{\Lambda(h_{i})+1}v_{\Lambda}\mid i=0,1,\ldots,l\}\subset
  V$$
  by \eqref{1.6}. According to Propositions 2.11 and 2.7, $F_{i}(-1)^{tk+1}
\mbox{{\bf 1}}\in U(\hat{\mathfrak{g}})x_{\theta}(-1)^{k+1}\mbox{{\bf
      1}}=N^{1}(k\Lambda_{0})$ for some $t\in \{1,2,3,4\}$. It
  then follows from \eqref{1.10}, \eqref{2.6} and  induction that
  $$Y_{M}^{\mu}(F_{i}(-1)^{tk+1}\mbox{{\bf
      1}},z)v_{\Lambda}\in V\big[\big[z^{1/r},z^{-1/r}\big]\big],$$
  and consequently
$$Y_{M}^{\mu}(F_{i}(-1)\mbox{{\bf
      1}},z)^{tk+1}v_{\Lambda}\in
  V\big[\big[z^{1/r},z^{-1/r}\big]\big]$$
  (since $[Y_{M}^{\mu}(F_{i}(-1)\mbox{{\bf
      1}},z_{1}),Y_{M}^{\mu}(F_{i}(-1)\mbox{{\bf
      1}},z_{2})]=0$ for $i=0,1,\ldots,l$). Hence
  \begin{equation}\label{2.24}
    f_{i}^{tk+1}v_{\Lambda}=\mbox{Res}_{z}\, z^{tk}Y_{M}^{\mu}(F_{i}(-1)
\mbox{{\bf 1}},z)^{tk+1}v_{\Lambda}\in V
  \end{equation}
  if $i\in \{1,\ldots,l\}$, and
  \begin{equation}\label{2.25}
    f_{0}^{tk+1}v_{\Lambda}=\mbox{Res}_{z}\, z^{\left(1-\frac{1}{r}\right)
(tk+1)-1}Y_{M}^{\mu}(F_{0}(-1)\mbox{{\bf
      1}},z)^{tk+1}v_{\Lambda}\in V,
  \end{equation}
  where $\mbox{$t=1,2,3$ or $4$}$ according to Proposition
  2.11. Since $[e_{i},f_{i}]=h_{i}$ and
  $h_{i}v_{\Lambda}=\Lambda(h_{i})v_{\Lambda}$, the representation
  theory of $\mathfrak{sl}(2,\mbox{{\bf C}})$ together with \eqref{2.6}, 
\eqref{2.24} and \eqref{2.25} imply that
  $$f_{i}^{\Lambda(h_{i})+1}v_{\Lambda}\subset
  V,\quad i=0,1,\ldots,l,$$
  as required.
  \end{proof}

\smallskip

It is easy to see that only some minor changes are needed in
order to make the above arguments also work in the case of the trivial
diagram automorphism $\mbox{id}_{\mathfrak{g}}$. In fact, for the
trivial twist most of these arguments can be simplified and one
obtains in this way the corresponding theorems in the 
untwisted case (cf.\@~\cite[Theorems 5.9 \& 5.14]{MP2}). Then one can 
extend the previous results to arbitrary finite-order
automorphisms, as shown by the following
\begin{theorem}
  Let $\sigma ,\psi \in \mbox{{\em Aut}}(\mathfrak{g})$ be such that 
  $\sigma \psi=\psi \mu \exp\!\big(\mbox{{\em ad}}
\!\left(\!\tfrac{2\pi
      i}{T}h\!\right)\!\big)$ for some appropriate
    $h\in \mathfrak{t}_{[0]}$ and some diagram automorphism $\mu$
    defined with respect to $(\mathfrak{t},\Pi)$, where $T$ is the order of $\sigma$. Let $M^{1}(\Lambda)$
    be the maximal
    submodule of the level $k$ Verma
    $\hat{\mathfrak{g}}[\sigma]$-module $M(\Lambda)$ with highest weight
  $\Lambda\in P_{+}$. Then
  \begin{equation}\label{2.26}
    \overline{\psi (R)}_{\sigma}M(\Lambda)=M^{1}(\Lambda),
  \end{equation} 
  where $R$ is as in \eqref{2.15} and $\overline{\psi
    (R)}_{\sigma}=\mathbf{C}\mbox{-{\em span}}\,\big\{\text{{\em Res}}_{z}\,
z^{n}Y_{M(\Lambda)}^{\sigma}(r,z)\mid r\in \psi (R),n\in \tfrac{1}{T}
\mathbf{Z}\big\}$.
\end{theorem}

\begin{proof}
  Note first that by Proposition 2.5 there exist indeed $\psi \in
  \mbox{Aut}(\mathfrak{g})$ and $h\in \mathfrak{t}_{[0]}$ as
  specified in the assumptions. Set
  $\mathfrak{t}'=\psi
  (\mathfrak{t})$, $\mu'=\psi \mu \psi^{-1}$, $h'=\psi (h)$, $R'=\psi (R)$, 
  and recall the automorphism $\bar{\mu}$ of order $r$ of
  the Dynkin diagram of $\mathfrak{g}$. Then 
  $\mathfrak{t}'$ is a $CSA$ of $\mathfrak{g}$ and 
  $\mu'$ is the diagram automorphism of $\mathfrak{g}$
  induced by the automorphism $\bar{\mu}$ with respect to
  $(\mathfrak{t}',\Pi)$. Let $\mathfrak{t}'=\coprod_{j\in
    \mathbf{Z}_{r}}\mathfrak{t}_{[j]}'$ be the $\mu'$-decomposition of
  $\mathfrak{t}'$ and notice that $h'\in \mathfrak{t}_{[0]}'$ and that
  $\sigma=\mu' \exp\!\big(\mbox{ad}\!\left(\!\tfrac{2\pi
      i}{T}h'\right)\!\big)$. Moreover, $R'$ is $\exp\!\big(\mbox{ad}\!
\left(\!\tfrac{2\pi
      i}{T}h'\right)\!\big)$-invariant and it obviously satisfies
    $\mathfrak{g}(n)R'=0$ for all integers $n\ge 1$. Then \eqref{2.26} 
    follows by combining Proposition 2.3 with Theorem 2.12 applied to the 
data $(\mathfrak{t}',\Pi ,R',\mu')$. 
\end{proof}

\smallskip

We now obtain the following characterization of
standard modules in terms of (actions of) irreducible loop $\tilde
{\mathfrak{g}}[\sigma]$-modules:
\begin{theorem}
  Let $\sigma,\psi \in \mbox{{\em Aut}}(\mathfrak{g})$ be as in Theorem 2.13, 
and let $M$ be a
  highest weight $\hat{\mathfrak{g}}[\sigma]$-module of level $k$. Then
  $M$ is a standard module if and only if $\overline{\psi (R)}_{\sigma}$
  annihilates $M$. 
\end{theorem}
\begin{proof}
  By Theorem 2.13 it suffices to consider only the case when
  $\sigma=\mu$ and $\psi=\mbox{id}_{\mathfrak{g}}$. Let $\Lambda$ be
  the highest weight of M and let $v_{\Lambda}\in M$ be a highest weight
  vector. If $\bar{R}_{\mu}M=0$, then one can argue as in the proof of
  Theorem 2.12 to show that for $0\le i\le l$ one has
  $f_{i}^{tk+1}v_{\Lambda}\in \bar{R}_{\mu}M=0$ for some $t\in
\{1,2,3,4\}$, so that $M$ is necessarily integrable. The converse is
given by Theorem 2.9. 
\end{proof}

\end{document}